\RequirePackage{ifpdf}
\ifpdf 
\documentclass[pdftex]{sigma}
\else
\documentclass{sigma}
\fi

\begin{document}

\allowdisplaybreaks

\renewcommand{\PaperNumber}{021}

\FirstPageHeading

\ShortArticleName{On the Existence of a Codimension 1 Completely Integrable Totally Geodesic Distribution}

\ArticleName{On the Existence of a Codimension 1 Completely\\ Integrable
Totally Geodesic Distribution\\ on a Pseudo-Riemannian Heisenberg Group}

\Author{Wafaa BATAT~$^\dag$ and Salima RAHMANI~$^{\ddag\S}$}

\AuthorNameForHeading{W.~Batat and S.~Rahmani}

\Address{$^\dag$~Ecole Normale Sup\'{e}rieure de L'Enseignement Technologique d'Oran,\\
\hphantom{$^\dag$}{}~D\'{e}partement de Math\'{e}matiques et Informatique, B.P. 1523 El M'Naouar Oran, Algeria}
\EmailD{\href{mailto:wafa.batat@enset-oran.dz}{wafa.batat@enset-oran.dz}, \href{mailto:wafa.batat@usc.es}{wafa.batat@usc.es}}

\Address{$^\ddag$~Laboratoire de Math\'ematiques-Informatique et Applications, Universit\'e de Haute Alsace,\\
\hphantom{$^\ddag$}{}~68093 Mulhouse Cedex, France}
\EmailD{\href{mailto:ramses161616@yahoo.fr}{ramses161616@yahoo.fr}}

\Address{$^\S$~Ecole Doctorale de Syst\`emes Dynamiques et G\'eom\'etrie,
D\'epartement de Math\'ematiques,\\ 
\hphantom{$^\S$}{}~Facult\'e des Sciences,
Universit\'e des Sciences et de la Technologie d'Oran,\\
\hphantom{$^\S$}{}~B.P. 1505  Oran El M'Naouer, Algeria}

\ArticleDates{Received December 23, 2009, in f\/inal form February 23, 2010;  Published online February 28, 2010}

\Abstract{In this note we prove that the Heisenberg group with a left-invariant pseudo-Riemannian metric admits a completely integrable totally geodesic distribution of codimension 1. This is on the
contrary to the Riemannian case, as it was proved by T.~Hangan.} 

\Keywords{Heisenberg group; pseudo-Riemannian
metrics; geodesics; codimension 1 distributions; completely integrable
distributions}

\Classification{58A30; 53C15; 53C30}

\section{Introduction}

Given an $n$-dimensional $C^{\infty}$ manifold $M$, a distribution $D$ of codimension $k$ of $M$ is an
$(n-k)$-dimensional subbundle of the tangent space $TM$.
A distribution $D$ of a pseudo-Riemannian manifold $\left( M,g\right) $ is
called \textit{totally geodesic} if every geodesic tangent to~$D$ at some
point remains everywhere tangent to~$D$. Among 2-step nilpotent Lie groups with left-invariant metrics, the Heisenberg group is of particular signif\/icance. In~\cite{R}, the second author proved that there are three nonisometric left-invariant Lorentzian metrics on the 3-dimensional Heisenberg group, one of which is f\/lat. This is a strong contrast to the Riemannian case in which there is only one (up to positive homothety) and it is not f\/lat.

Another major dif\/ference to the Riemannian case appears. Indeed, T. Hangan
proved that the Heisenberg group endowed with a left-invariant Riemannian
metric does not admit any codimension 1 completely integrable totally geodesic distribution
(see~\cite{Han1} for the three-dimensional case and~\cite{Han2} for the
higher dimensional case). Recently, in~\cite{RR} the second author and N.~Rahmani proved that the 3-dimensional Heisenberg group admits a left-invariant Lorentzian metric for which there exists a codimension 1 completely integrable totally geodesic distribution.
It is an interesting problem to investigate whether and
to what extent results valid in Riemannian Heisenberg group can be extended to the pseudo-Riemannian
case. In fact, in a high-dimensional Heisenberg group there are dif\/ferent types of pseudo-Riemannian metrics while there is still essentially only one (up to positive homothety) Riemannian metric. So,
the purpose of this paper is to prove, contrary to the Riemannian case, the
existence of a codimension 1 completely integrable totally geodesic
distribution on a pseudo-Riemannian Heisenberg group of
dimension $2p+1$.

\section{Geodesics and totally geodesic distributions\\ on Heisenberg group}

Consider $\mathbb{R}^{2p+1}$ with the elements of the form $X=(x_{1},\dots,x_{p},y_{1},\dots,y_{p},z)$.
Def\/ine the product on $%
\mathbb{R}^{2p+1},$ by%
\[
X\tilde{X}=\left(x_{1}+\tilde{x}_{1},\dots,x_{p}+\tilde{x}_{p},y_{1}+\tilde{y}_{1},\dots,y_{p}+%
\bar{y}_{p},z+\tilde{z}-\overset{p}{\underset{i=1}{\sum }}x_{i}\tilde{y}%
_{i}\right),
\]
where $X=(x_{1},\dots,x_{p},y_{1},\dots,y_{p},z)$, $\tilde{X}=(\tilde{x}_{1},\dots,\tilde{x}%
_{p},\tilde{y}_{1},\dots,\tilde{y}_{p},\tilde{z}).$

Let $H_{2p+1}=\left(\mathbb{R}^{2p+1},g\right) $ be the $(2p+1)$-dimensional Heisenberg group with the left-invariant
pseudo-Riemannian metric $g$ def\/ined by%
\begin{equation}
g=-\overset{p}{\underset{i=1}{\sum }}\left( dx_{i}\right) ^{2}+\overset{p}{%
\underset{i=1}{\sum }}\left( dy_{i}\right) ^{2}+\left( dz+\overset{p}{%
\underset{i=1}{\sum }}x_{i}dy_{i}\right) ^{2}.  \label{g}
\end{equation}%

With respect to the basis of coordinate vector f\/ields $\big\{ \partial
_{x_{i}}=\frac{\partial }{\partial x_{i}},\partial _{y_{i}}=\frac{\partial }{%
\partial y_{i}},\partial _{z}=\frac{\partial }{\partial z}\big\} $, $i=1,2,\dots,p$ for which \eqref{g} holds, the nonvanishing metric components are%
\begin{gather}
g_{_{\partial _{x_{i}}\partial _{x_{i}}}}=-1,\qquad g_{_{\partial
_{y_{i}}\partial _{y_{i}}}}=1+\left( x_{i}\right) ^{2},\qquad
g_{_{\partial _{y_{i}}\partial _{z}}}=x_{i}\qquad \mathrm{for}\quad 1\leq i\leq p,
\nonumber\\
g_{_{\partial _{z}\partial _{z}}}=1,\qquad   g_{_{\partial _{y_{i}}\partial
_{y_{j}}}}=x_{i}x_{j}\qquad \mathrm{for}\quad 1\leq i\neq j\leq p.
\label{gij}
\end{gather}

The geodesics equations are obtained in \cite{Eb} (for the Riemannian case)
and in~\cite{ParCor} and~\cite{Gued} (for the pseudo-Riemannian case), in
more general context of two-step nilpotent Lie groups with a~left-invariant
metric. Jang and Parker~\cite{Jang} later gave explicit formulas for geodesics
on the 3-dimensional Heisenberg group endowed with a left-invariant Lorentzian metric.
A corresponding description was made in~\cite{Han1} in the Riemannian case.

Here we will f\/ind the geodesics equations of $H_{2p+1}$ for the left-invariant
pseudo-Riemannian metric~\eqref{g} by integrating their Euler--Lagrange
equations. Notice that Hangan  obtained in~\cite{Han2} all
geodesics of the Riemannian Heisenberg group of dimension $2p+1$.

So, let $\gamma (t)$ be a geodesic curve on $H_{2p+1},$ determined,
with respect to the basis of coordinate vector f\/ields $\big\{ \frac{%
\partial }{\partial x_{1}},\dots,\frac{\partial }{\partial x_{p}},\frac{%
\partial }{\partial y_{1}},\dots,\frac{\partial }{\partial y_{p}},\frac{%
\partial }{\partial z}\big\} $, by its components $(x(t),y(t),z(t))$ where
$x(t)=(x_{1}(t),\dots,x_{p}(t))$, $y(t)=(y_{1}(t),\dots,y_{p}(t))$ and $z(t)\in
\mathbb{R}.$ Using~\eqref{gij}, the Lagrangian associated to the metric~\eqref{g} is
\[
L=\frac{1}{2}\left\{ -\overset{p}{\underset{i=1}{\sum }}\left( x_{i}^{\prime
}\right) ^{2}+\overset{p}{\underset{i=1}{\sum }}\left( y_{i}^{\prime
}\right) ^{2}+\left[ z^{\prime }+\overset{p}{\underset{i=1}{\sum }}%
x_{i}y_{i}^{\prime }\right] ^{2}\right\} ,
\]%
with $x^{\prime }(t)=(x_{1}^{\prime }(t),\dots,x_{p}^{\prime }(t))$ and $%
y^{\prime }(t)=(y_{1}^{\prime }(t),\dots,y_{p}^{\prime }(t))$. Hence, the
corresponding Euler--Lagrange equations are
\begin{gather}
x^{\prime \prime }+\left( z^{\prime }+ \sum_{i=1}^p
x_{i}y_{i}^{\prime }\right) y^{\prime }=0, \qquad
y^{\prime \prime }+\left( z^{\prime }+ \sum_{i=1}^p
x_{i}y_{i}^{\prime }\right) x^{\prime }=0, \nonumber\\
z^{\prime \prime }+ \sum_{i=1}^p \big( x_{i}^{\prime
}y_{i}^{\prime }+x_{i}y_{i}^{^{\prime \prime }}\big) =0,
  \label{sys}
\end{gather}%
and we will integrate the system \eqref{sys} with the initial conditions%
\begin{alignat*}{4}
& x(0)=x_{0}, \qquad && y(0)=y_{0}, \qquad && z(0)=z_{0}, & \\
& x^{\prime }(0)=x_{0}^{\prime },\qquad && y^{\prime }(0)=y_{0}^{\prime },\qquad && z^{\prime }(0)=z_{0}^{\prime }.&
\end{alignat*}
The f\/irst integration of the last equation of the Euler--Lagrange
equations, reduces the system~\eqref{sys}~to
\begin{gather}
x^{\prime \prime }+\alpha y^{\prime }=0,\qquad
y^{\prime \prime }+\alpha x^{\prime }=0, \qquad
z^{\prime }+ \sum_{i=1}^p  x_{i}y_{i}^{\prime }=\alpha ,
\label{sysred}
\end{gather}
where $\alpha =z_{0}^{\prime }+ \sum\limits_{i=1}^p
x_{i}(0)y_{i}^{\prime }(0)$ is a real constant. For the system~\eqref{sysred}
we distinguish two cases.

If $\alpha =0,$ the geodesics curves are given by
\begin{gather}
x(t)=x_{0}^{\prime }t+x_{0}, \qquad
y(t)=y_{0}^{\prime }t+y_{0}, \nonumber\\
z(t)=z_{0}- \sum_{i=1}^p \left( x_{i}^{\prime }(0)\frac{t}{2}+x_{i}(0)\right) y_{i}^{\prime }(0)t.
 \label{geod0}
\end{gather}

If $\alpha \neq 0,$ the geodesics curves are given by%
\begin{gather}
x(t)=\frac{1}{\alpha }\big( \sinh(\alpha t)x_{0}^{\prime }-\cosh(\alpha
t)y_{0}^{\prime }+\alpha x_{0}+y_{0}^{\prime }\big) , \nonumber\\
y(t)=\frac{1}{\alpha }\big( -\cosh(\alpha t)x_{0}^{\prime }+\sinh(\alpha
t)y_{0}^{\prime }+\alpha y_{0}+x_{0}^{\prime }\big) ,\nonumber \\
z(t)=\frac{1}{\alpha ^{2}}\overset{p}{\underset{i=1}{\sum }}\left( \alpha
x_{i}(0)+y_{i}^{\prime }(0)\right) \left( x_{i}^{\prime }(0)\cosh(\alpha
t)-y_{i}^{\prime }(0)\sinh(\alpha t)\right) \nonumber\\
\phantom{z(t)=}{}
-\frac{1}{2\alpha ^{2}}
\sum_{i=1}^p x_{i}^{\prime }(0)y_{i}^{\prime }(0)\cosh(2\alpha t) \nonumber\\
\phantom{z(t)=}{}+
\frac{1}{4\alpha ^{2}}\sum_{i=1}^p \big[ (
x_{i}^{\prime }(0)) ^{2}+( y_{i}^{\prime }(0)) ^{2}\big]
\sinh(2\alpha t)+\left( \alpha -\frac{1}{2\alpha }\sum_{i=1}^p \big[ ( x_{i}^{\prime }(0)) ^{2}-( y_{i}^{\prime
}(0)) ^{2}\big] \right) t  \nonumber\\
\phantom{z(t)=}{}+
z_{0}-\frac{1}{2\alpha ^{2}}\sum_{i=1}^p x_{i}^{\prime
}(0)\left( 2\alpha x_{i}(0)+y_{i}^{\prime }(0)\right) .%
  \label{geod}
\end{gather}

It must be noted that the pseudo-Riemannian metric $g$ on a $(2p+1)$-dimensional Heisenberg group
is geodesically complete. This is a general truth for pseudo-Riemannian 2-step nilpotent Lie groups. In fact,
in \cite{Gued1} Guediri proved that every left-invariant pseudo-Riemannian metric on a 2-step nilpotent Lie groups
is geodesically complete.

We now prove the following result

\begin{theorem}
The pseudo-Riemannian metric \eqref{g} admits a codimension $1$ completely
integrable totally geodesic distribution on the Heisenberg group
$H_{2p+1}.$
\end{theorem}

\begin{proof} Let $\gamma $ be a geodesic curve on $%
H_{2p+1}$ with $\gamma (0)=\left( x_{0},y_{0},z_{0}\right) $ and $%
\gamma ^{\prime }(0)=\left( x_{0}^{\prime },y_{0}^{\prime },z_{0}^{\prime
}\right) $. Consider in $H_{2p+1}$ the 1-form $\theta =\overset{p}{%
\underset{i=1}{\sum }}\left( dx_{i}-dy_{i}\right) .$

If $z_{0}^{\prime }+\sum\limits_{i=1}^p x_{i}(0)y_{i}^{\prime
}(0)=0$ then $\gamma $ is given by \eqref{geod0} and we have%
\begin{gather*}
x^{\prime }(t)  = x_{0}^{\prime }, \qquad
y^{\prime }(t)  = y_{0}^{\prime}.
\end{gather*}

If $z_{0}^{\prime }+\sum\limits_{i=1}^p x_{i}(0)y_{i}^{\prime
}(0)\neq 0$ then $\gamma $ is given by \eqref{geod} and we have%
\begin{gather*}
x_{i}^{\prime }(t)  = \left( \cosh(\alpha t)x_{i}^{\prime }(0)-\sinh(\alpha
t)y_{i}^{\prime }(0)\right) , \qquad
y_{i}^{\prime }(t)  = \left( -\sinh(\alpha t)x_{i}^{\prime }(0)+\cosh(\alpha
t)y_{i}^{\prime }(0)\right) ,
\end{gather*}%
for all $i\in \left\{ 1,2,\dots,p\right\} .$ From \eqref{geod0} and \eqref{geod}
one gets
\begin{gather*}
\theta \left( \gamma ^{\prime }(t)\right) =\sum_{i=1}^p
\left( x_{i}^{\prime }(0)-y_{i}^{\prime }(0)\right) \qquad \text{if} \quad
z_{0}^{\prime }+\sum_{i=1}^p x_{i}(0)y_{i}^{\prime}(0)=0 , \\
\theta \left( \gamma ^{\prime }(t)\right) =\left[ \cosh(\alpha t)+\sinh(\alpha t)%
\right] \left( \sum_{i=1}^p \left( x_{i}^{\prime
}(0)-y_{i}^{\prime }(0)\right) \right) \qquad \text{if} \quad z_{0}^{\prime }+\sum_{i=1}^p x_{i}(0)y_{i}^{\prime }(0)\neq 0,
\end{gather*}
for all $t\in\mathbb{R}.$ This means that, in both cases, if $\theta \left( \gamma ^{\prime
}(0)\right) =0$ then $\theta \left( \gamma ^{\prime }(t)\right) =0$ for all $%
t\in\mathbb{R}.$ In other words, the kernel of the given 1-form $\theta $ def\/ines a
codimension 1 totally geodesic distribution. Since $\theta$ is closed, the distribution $D=\ker \theta $ is
involutive and from the Frobenius theorem it follows that it is completely integrable.
\end{proof}

Notice that since $\theta $ is the dif\/ferential of the function $f =\sum\limits_{i=1}^p\left( x_{i}-y_{i}\right)$, the submanifolds given by $f=\textrm{const}$ are the
maximal integral submanifolds of $D$.

\medskip

\noindent{\bf Remark.}
Since the dimension of the center of a Heisenberg group is~1, it seems natural to pose the study of
the existence of totally geodesic distributions on 2-step nilpotent Lie groups endowed with a pseudo-Riemannian
metric with index equal to the dimension of the center of the Lie algebra.

\subsection*{Acknowledgements}

The authors wish to express their gratitude toward the referees for their valuable remarks and also would like to thank Professor Jos\'{e} Antonio Oubi\~na for his fruitful comments on the revised version of this paper. The f\/irst author was supported by ENSET d'Oran and the second author is indebted to Doctoral School of Dynamical Systems and Geometry (U.S.T.~Oran) for its f\/inancial support during the elaboration of this work.

\pdfbookmark[1]{References}{ref}
\LastPageEnding

\end{document}